\documentclass[11pt,oneside,a4paper]{article} 
\usepackage{lscape}
\usepackage{color}
\usepackage{amsfonts}
\usepackage{latexsym}
\usepackage{graphicx}
\usepackage{url} 
\usepackage[english]{babel}
\usepackage{babelbib}
\usepackage{algorithm, algpseudocode}

\usepackage{times}             
\usepackage[latin1, applemac]{inputenc} 
\usepackage{graphicx}          

\usepackage{textcomp}
\usepackage{mathptmx}

\usepackage{amsmath, amsbsy, amssymb, amsthm}
\newtheorem*{myLemma}{Lemma}
\newtheorem*{myProof}{Proof}
\begin{document}

\title{Factor PD-Clustering }
\author{{Mireille Gettler Summa$^1$, Francesco Palumbo$^2$, Cristina Tortora$^{1,2}$.  }\\
{$^1$Université Paris Dauphine CEREMADE CNRS}\\
{$^2$Università di Napoli Federico II }}
\date{}
\maketitle

\section*{Abstract}
Factorial clustering methods have been developed in recent years thanks to the improving of computational power. These methods perform a linear transformation of data and a clustering on transformed data optimizing a common criterion. Factor PD-clustering is based on Probabilistic Distance clustering (PD-clustering).
 PD-clustering is an iterative, distribution free, probabilistic, clustering method.
Factor PD-clustering makes a linear transformation of original variables into a reduced number of orthogonal ones using  a common criterion with PD-Clustering. This paper demonstrates that Tucker3 decomposition permits to obtain this transformation. Factor PD-clustering exploits alternatively  Tucker3 decomposition and PD-clustering on transformed data until convergence is achieved. 
This method can significantly improve the algorithm performance; large datasets can thus be partitioned into clusters with increasing  stability and robustness of the results.

\section*{keyword}
Multivariate analysis, Exploratory data analysis, Clustering, Factorial clustering, Non hierarchical iterative clustering.




\section{Introduction}
\label{intro}
In a wide definition Cluster Analysis is a multivariate analysis technique that seeks to organize information about
variables in order to discover  homogeneous groups, or ``clusters'' into data. 
The presence of groups in data depends on the association structure over the data.
Clustering algorithms aim at finding homogeneous groups with respect to their
association structure among variables. Proximity measures or distances can be properly used to
separate homogeneous groups.  \\
A measure of the homogeneity of a group is the variance.
Dealing with numerical linearly
independent  variables, clustering problem consists in minimising the sum of the squared Euclidean 
distances within classes: within groups deviance.

{\it ``The term cluster analysis refers to an entire process where clustering maybe only a step"} \cite{Gor99}.
According to Gordon's definition cluster analysis can be sketched in three  main stages:
\begin{itemize}
\item transformation of data into a similarity/dissimilarity matrix;
\item clustering;
\item validation.
\end{itemize}
Transformation of data into similarity/dissimilarity measures depends on data type. On the transformed matrix a clustering method can be applied. Clustering methods can be divided into three main types:  hierarchical, non hierarchical and fuzzy \cite{WedKam99}. Non hierarchical clustering methods are considered in this paper. Among them the most well-known and used method is k-means. It is an iterative method that starts with a random initial partition of units and keeps reassigning the units into clusters based on the squared distances between the units and the clusters' centers until the convergence is reached. Interested readers can refer to \cite{Gor99}. Major k-means issues are that clusters can be sensitive to the choice of the initial centers and that the algorithm could converge to local minima. 

 The choice of the number of clusters is a well known problem of non hierarchical methods, this problem will not be dealt with in this paper where the number of clusters is assumed as a priori  known. \\
Non hierarchical clustering methods performance can be strongly affected by the dimensionality. Let us consider  an $n \times J$ data matrix $X$, with $n$ number of units and $J$ number of variables. Non hierarchical methods easily deal with large $n$, however they can fail when $J$ becomes large or very large and when the variables are correlated. They do not converge or they converge into a different solution at each iteration.
To cope with these issues, French school of  data analysis \cite{LebMor84} suggested a strategy to improve the overall quality of clustering that consists in two phases: variables transformation through a factorial method and clustering method on transformed variables.
Arabie and Hubert in 1996 \cite{AraHub96} fourthly formalized the method  and called it {\it  tandem analysis}.\\
The choice of the factorial method is an important and tricky phase because it will affect the results. 
Principal factorial methods are \cite{LerRou04}:

\begin{itemize}
\item quantitative data;
\begin{itemize}
\item  Principal Component Analysis (PCA);
\end{itemize}
\item binary data;
\begin{itemize}
\item  Principal Component Analysis (PCA);
\item Correspondence Analysis (CA);
\item  Multiple Correspondence Analysis (MCA);
\end{itemize}
\item nominal data;
\begin{itemize}
\item   Multiple Correspondence Analysis (MCA);
\end{itemize}
\end{itemize}

The second phase of the {\it  tandem analysis} consists in applying clustering methods.\\
Tandem analysis exploits the factor analysis capabilities that consist in obtaining a reduced number of uncorrelated variables which are linear transformation of the original ones. This method gives more stability to the results and makes the procedure faster. However tandem analysis minimises two different functions that can be in contrast and  the first factorial step can in part obscure or mask the clustering structure.  \\
This technique has the advantage of working with a reduced number of variables that are orthogonal and ordered with respect to the borrowed information. Moreover dimensionality reduction permits to visualize the cluster structure in two or three dimensional factorial space \cite{PalVis08}.
To cope with these issues Vichi and Kiers \cite{VicKie01} proposed Factorial k-means analysis for two-way data. The aim of this method is to identify the best partition of the objects and to find a subset of factors that best describe the classification according to the least squares criterion. Two steps Alternating Least Squares algorithm (ALS) based on  solves this problem.
The advantage of Factorial k-means is that the two steps optimize a single objective function. 
However the k-means  algorithm itself, and as a consequence the tandem analysis and Factorial k-means, is based on the arithmetic mean that gives rise to unsatisfactory solutions when clusters have not spherical shape.

Probabilistic clustering methods may allow us to obtain better results under this condition because they assign a statistical unit to a cluster according to a probability function that can be independently defined with respect to the arithmetic mean. \\
Probabilistic Distance clustering (PD-clustering) \cite{BenIyi08} is an iterative, distribution free, probabilistic, clustering method. PD-clustering assigns units to a cluster according to their probability of belonging to the cluster, under  the constraint  that the product between the probability and the distance of each point to any cluster center is a constant.\\
When the number of variables is large and variables are correlated, PD-clustering becomes unstable and the correlation between variables can hide the real number of clusters. A linear transformation of original variables into a reduced number of orthogonal ones using  common criteria with PD-clustering can significantly improve the algorithm performance. 
The objective of this paper is to introduce an improved version of PD-clustering called Factor PD-clustering (FPDC).

The paper has the following structure:
section 2: detailed presentation of PD-clustering method;
section 3: presentation of our suggestion for a Factor PD-clustering method;
section 4: application of Factor PD-clustering on a simulated case study and comparison with k-means.

\section{Probabilistic Distance Clustering}\label{PDclustering}

PD-clustering  is a non hierarchical algorithm that assigns units to clusters according to their belonging probability to the cluster. According to Ben-Israel and Iyigun \cite{BenIyi08} notation we introduce PD-clustering. Given some random centers, the probability of any point to belong to each class is assumed to be inversely proportional to the distance from the centers of the clusters.
Given an $X$  data matrix with $n$ units and $J$ variables, given $K$ clusters that are assumed not empty,
 PD-Clustering is based on two quantities:  the distance of each data point $x_i$ from the $K$ cluster centers $c_k$, $d(x_i,c_k)$, and the probabilities for each point to belong to a cluster, $p(x_i,c_k)$ with $k=1,\ldots , K$ and $i=1,\ldots, n$.
 The relation between  them is the basic assumption of the method. Let us consider the general term
 $x_{ij}$ of $X$ and a center matrix $C$, of elements $c_{kj}$ with $k=1,\ldots , K$, $i=1,\ldots, n$ and $j=1,\ldots, J$, their distance can be computed according to different criteria,  the squared norm is one of the most commonly used. The generic distance $d(x_i,c_k)$ represents the distance of the generic point $i$ to the generic center $k$. The probability $p(x_i,c_k)$ of each point to belong to a cluster can be computed according to the following assumption: the product between the distances and the probabilities is a constant depending on $x_i$: $F(x_i)$. \\
For short we use $p_{ik}=p(x_i,c_k)$ and $d_k(x_i)=d(x_i,c_k)$; PD-clustering basic assumption is expressed as:
\begin{equation}
\label{basicDC}{
p_{ik}d_k(x_i)=F(x_i).
}
\end{equation}
for a given value of  $x_i$ and for all $k=1, \ldots, K$.

At the decreasing of the point closeness from the cluster center the belonging probability of the point to the cluster decreases.
The constant depends only on the point and does not depend on the cluster $k$.\\
Starting from the \ref{basicDC} it is possible to compute $p_{ik}$:
\begin{eqnarray}
\label{pDC}{
p_{im}d_m(x_i)=p_{ik}d_k(x_i);\
p_{im}=\frac{p_{ik}d_k(x_i)}{d_m(x_i)};\
\forall m=1, \ldots, K
}
\end{eqnarray}
The term $p_{ik}$ is a probability so, under the constraint $\sum_{m=1}^K p_{im}=1$, the sum over $m$ of \ref{pDC} becomes:
\begin{eqnarray}
p_{ik}\sum_{m=1}^K \left(\frac{d_k(x_i)}{d_m(x_i)}\right)=1,\nonumber
\end{eqnarray}
\begin{eqnarray}
\label{pDC2}
p_{ik}=\left(\sum_{m=1}^K \left(\frac{d_k(x_i)}{d_m(x_i)}\right)\right)^{-1}=\frac{\prod_{m\neq k}d_m(x_i)}{\sum_{m=1}^K\prod_{k\neq m}d_k(x_i)}, k=1,\ldots K.
\end{eqnarray}
Starting from the \ref{basicDC} and using \ref{pDC2}
it is possible to define the value of the constant $F(x_i)$:
\begin{eqnarray}
F(x_i)=p_{ik}d_k(x_i), k=1,\ldots K, \nonumber
\end{eqnarray}
\begin{eqnarray}\label{eq:JDF}
F(x_i)=\frac{\prod_{m=1}^Kd_m(x_i)}{\sum_{m=1}^K\prod_{k\neq m}d_k(x_i)}.
\end{eqnarray}
The quantity $F(x_i)$, also called {\it Joint Distance Function} (JDF), is a measure of the closeness of $x_i$ from all clusters' centers. The JDF measures the classificability of the point $x_i$ with respect to the centers $c_k$ with $k=1, \ldots, K$. If it is equal to zero, the point coincides with one of the clusters' centers, in this case the point belongs to the class with probability $1$. If all the distances between the point $x_i$ and the centers of the classes are equal to $d_i$, $F(x_i)=d_i/k$ and all the belonging probabilities to each class are equal: $p_{ik}=1/K$. The smaller the JDF value, the higher the probability for the point to belong to one cluster. \\
The whole clustering problem consists in the identification of the centers that minimises the JDF. \\
Without loss of generality the PD-Clustering optimality criterium can be demonstrated according to $k=2$.
\begin{eqnarray}
\label{optimization}
&\min \left( d_1(x_i)p_{i1}^2 +d_2(x_i)p_{i2}^2 \right) \\
&\mbox{s.t.} & p_{i1}+p_{i2}=1 \nonumber \\
& & p_{i1}, p_{i2}\geq 0 \nonumber 
\end{eqnarray}
The probabilities are squared because it is a smoothed version of the original function.
The Lagrangian of this problem is:
\begin{eqnarray}
\mathcal{L}(p_{i1},p_{i2},\lambda)=d_1(x_i)p_{i1}^2+d_2(x_i)p_{i2}^2-\lambda(p_{i1}+p_{i2}-1)
\end{eqnarray}
Setting to zero the partial derivates with respect to $p_{i1}$ and $p_{i2}$, substituting the probabilities \ref{pDC2} and considering the principle $p_{i1}d_1(x_i)=p_{i2}d_2(x_i)$ we obtain the optimal value of the Lagrangian.
\begin{eqnarray}
\mathcal{L}(p_{i1},p_{i2},\lambda)=\frac{d_1(x_i)d_2(x_i)}{d_1(x_i)+d_2(x_i)}.
\end{eqnarray}
This value coincides with the JDF, the matrix of centers that minimises this principle minimises the JDF too.
Substituting the generic value $d_k(x_i)$ with $\left\|x_i-c_k\right\|$, we can find the equations of the centers that minimise the JDF (and maximize the probability of each point to belong to only one cluster).
\begin{eqnarray}\label{centri}
c_k=\sum_{i=1,\ldots,N}\left(\frac{u_k(x_i)}{\sum_{j=1,\ldots,N}u_k(x_j)}\right)x_i,
\end{eqnarray}
where
\begin{eqnarray}
\label{uk}
u_k(x_i)=\frac{p_{ik}^2}{d_k({x_i})}.
\end{eqnarray}
As showed before, the value of JDF at all centers $k$ is equal to zero and it is necessarily positive elsewhere. So the centers are the global minimiser of the JDF. Other stationary points may exist because the function is not convex neither quasi-convex, but they are saddle points.

There are alternative ways for modeling the relation between probabilities and distances, for example the probabilities can decay exponentially as distances increase. In this case the probabilities $p_{ik}$ and the distances $d_k(x_i)$ are related by:
\begin{eqnarray}
p_{ik} e^{d_k(x_i)} = E(x_i),
\end{eqnarray}
where $E(x_i)$ is a constant depending on $x_i$.\\
Many results of the previous case can be extended to this case by replacing the distance $d_k(x_i)$ with $e^{d_k(x_i)}$.
Interested readers are referred to Ben-Israel and Iyigun \cite{BenIyi08}.

The optimization problem presented in \ref{optimization} is the original version proposed by Ben-Israel and Iyigun.
Notice that in the optimization problem the probabilities $p_k$ are considered in squared form. 
The Authors affirm that it is possible to consider $d_k$ as well $d_k^2$. Both choices have some advantages and drawbacks.
Squared distances offer analytical advantages due to linear derivates.
Using simple distances endures more robust results and the optimization problem can be reconducted to a  Fermat-Weber location problem. The Fermat-Weber location problem aims at finding a point that minimises the sum of the Euclidean distances from a set of given points. This problem can be solved  with the Weiszfeld method \cite{Wei37}. Convergence of this method was established by modifying the gradient so that it is always defined \cite{Ku73}. The modification is not carried out in practice. The global solution is guaranteed only in case of one cluster. Dealing with more than one cluster, in practice, the method converges only for a limited number of centers depending on the data.

In this paper we consider the squared form:
\begin{eqnarray}\label{distanze}
d_k(x_i) =\sum_{j=1}^J(x_{ij} - c_{kj})^2 \label{distance1},  
\end{eqnarray}
where $k=1,\ldots,K$ and $i =1,\ldots , N$. Starting from the \ref{distanze} the distance matrix $D$ of order $n \times K$ is defined, where the general element is $d_k(x_i)$. The final solution $\hat{JDF}$ is obtained minimising the quantity:

\begin{eqnarray}
JDF = \sum_{i=1}^n  \sum_{k=1}^K d_k(x_i) p_{ik}^2 = \sum_{i=1}^n\sum_{j=1}^J\sum_{k=1}^K (x_{ij} - c_{kj})^2 p_{ik}^2, \label{distance2}
\end{eqnarray}
\begin{eqnarray}\label{JDFhat}
\hat{JDF} = \arg\min_{C;  P}{ \sum_{i=1}^n\sum_{j=1}^J\sum_{k=1}^K (x_{ij} - c_{kj})^2 p_{ik}^2}.
\end{eqnarray}
Where $c_k$ is the generic center and $d_k(x_i)$ is defined in \ref{distance1}.

\vspace{\baselineskip}

The solution of PD-clustering problem can be obtained through an iterative algorithm.\\

\begin{algorithm}\begin{footnotesize}
\caption{Probabilistic Distance Clustering Function}\label{PDC:algo}
\begin{algorithmic}[1]
\Function{PDC}{$X,K$}
\State $C\gets \mathrm{rand}(K,J)$ \Comment{Matrix $C_{K,J}$ is randomly initialised}
\State $JDF \gets 1/\mathrm{eps}$           \Comment{$JDF$ is initialised to the maximum}
\State $D\gets 0$                             \Comment{Initialise the array $D$, of dimension $n\times K$, to $0$}
\State $p \gets \frac{1}{k}$                             \Comment{Initialise to $\frac{1}{k}$ the probability vector $p$ of $n$ elements} 
\Repeat
\For{ $k=1,K$}
\State $D_k \gets distance(X,C(k))$ \Comment{$D_k$ distances of all units from the centre $k$ according to formula \ref{distanze}}
\EndFor

\State $JDF0 \gets JDF$                \Comment{Current $JDF$ is stored in $JDF0$} 
\State $C \gets C^{*}$                     \Comment{Centres are updated according to formula \ref{centri}}

\State $JDF\gets jdf(D)$              \Comment{$\gets jdf(D)$ implements the formula \ref{eq:JDF}}
\Until{$JDF0 > JDF$}
\State $P\gets compp(D)$        \Comment{function $compp$ implements the formula \ref{pDC2}}

\Return{$C, P, JDF$}
\EndFunction

\end{algorithmic}\end{footnotesize}
\end{algorithm}

The algorithm convergence is demonstrated in \cite{Iyi07}.\\
Each unit is then assigned to the $k^{th}$ cluster according to the highest probability that is computed a posteriori using the formula in
 equation \ref{pDC2}.

\section{Factor PD-Clustering}
When the number of variables is large and variables are correlated, PD-Clustering becomes very unstable and the correlation between variables can hide the real number of clusters. A linear transformation of original variables into a reduced number of orthogonal ones can significantly improve the algorithm performance. Combination of PD-Clustering and variables linear transformation implies a common criterion. 

This section shows how the Tucker3 method \cite{Kro08} can be properly adopted for the  transformation into the Factor PD-Clustering; an algorithm is then proposed to perform the method.

\subsection{Theoretical approach to Factor PD-clustering}\label{sec:3.1}

Firstly we demonstrate that the minimization problem in \ref{distance2} corresponds to the Tucker3 decomposition of the 
distance matrix  $G$ of general elements $g_{ijk}=|x_{ij}-c_{kj}|$. 
It is a 3-way matrix $n \times J \times K$  where  $n$ is the number of units, $J$ the number of variables and $K$ the occasions.
For any $c_k$ with $k=1, \ldots, K$, a $G_k$ $n \times J$ distances matrix is defined. 
In matrix notation:
\begin{eqnarray}\label{dist}
G_k=X-hc_k
\end{eqnarray}
where $h$ is an $n\times 1$ column vector with all terms equal to 1; $X$ and $c_k$ ($k=1,\ldots, K$) have been already defined
in section \ref{PDclustering}.

Tucker3  method decomposes the matrix $G$ in three components, one for each mode, in a full core array $\Lambda$ and in an error term $E$. 
\begin{eqnarray}\label{eq:g}
g_{ijk}=\sum_{r=1}^R\sum_{q=1}^Q \sum_{s=1}^S \lambda_{rqs}(u_{ir}b_{jq}v_{ks})+ e_{ijk},  
\end{eqnarray}
where $\lambda_{rqs}$ and  $e_{ijk}$ are respectively  
the general terms of the three way matrix $\Lambda$ of order  $R\times S\times Q$ and $E$ of order $n\times J \times K$;

 $u_{ir}$, $b_{jq}$ and $v_{ks}$  are respectively  
the general terms of the matrix $U$ of order $n\times R$, $B$ of order $J\times Q$ and $V$ of order $K\times S$, with $ i=1,\ldots, n$, $j=1,\ldots, J$, $k=1,\ldots, K$.

As in all factorial methods, factorial axes in Tucker3 model are sorted according to explained variability. The first factorial axes explain the greatest part of the variability, latest factors are influenced by anomalous data or represent the ground noise. For this reason the choice of a number of factors  lower than the number of variables makes the method externally robust.
According to \cite{KieKin03} the choice of the parameters $R$, $Q$ and $S$ is a ticklish
problem because they define the overall explained variability. The interested readers 
are referred to \cite{Kro08} for the theoretical aspects concerning this choice. 
We use an heuristic approach to cope with this crucial issue: we choose the minimum number of factors that corresponds to a significant value of the explained variability.

\noindent
The coordinates $x^{*}_{iq}$ of the generic unit $x_i$ 
into the space of variables obtained through Tucker3 decomposition 
are obtained by the following expression:
\begin{eqnarray}
\label{coordinate}
x^{*}_{iq}=\sum_{j=1}^J x_{ij}b_{jq}.
\end{eqnarray}

\noindent
Finally on these $x_{iq}^{*}$ coordinates a PD-Clustering is applied in order to solve the clustering problem. 

Let us start considering the expression \ref{distance2};
it is worth noting that minimising the quantity:
\begin{eqnarray}\label{obb}
\textstyle{
JDF= \sum_{i=1}^n\sum_{j=1}^J\sum_{k=1}^K(x_{ij}-c_{kj})^2 p_{ik}^2}\hspace{0.1\textwidth} \mathrm{s.t. }\quad \sum_{i=1}^n\sum_{k=1}^K{p_{ik}^2} \leq n,
\end{eqnarray}
is equivalent to  compute the maximum of $-\sum_{i=1}^n\sum_{j=1}^J\sum_{k=1}^K(x_{ij}-c_{kj})^2 p_{ik}^2$,
under the same constraints.

Taking into account the Proposition 1 (proof in \ref{App:1}) and the following lemma, we demonstrate that the Tucker3 decomposition 
is a consistent linear variable transformation that determines the best subspace according to the PD-clustering criterion.\\
{\bf Proposition 1} {\it Given an unknown matrix $B$ of generic element $b_{im}$ and a set of  coefficients $0 \leq  \psi_{im} \leq 1$, with $m=1,\ldots, M$ and $i=1, \ldots, n$.
Maximising}
$$-\sum_{m=1}^{M}\sum_{i=1}^nb_{im}\psi_{im}^2,$$ $\mathrm{s.t. } \sum_{m=1}^M\sum_{i=1}^n \psi^2_{im} \leq n$ 
{\it is equivalent to solve the equation}
$$\sum_{m=1}^{M}\sum_{i=1}^nb_{im}\psi_{im }= \mu\sum_{m=1}^M\sum_{i=1}^n\psi_{im},$$ 
{\it where $\mu \geq 0 $.\\}

\begin{myLemma} Tucker3 decomposition permits to define the best subspace for the PD-clustering.
\end{myLemma}

We consider the proposition of the {\em Proposition 1} where:
\begin{eqnarray}
M&=&K \nonumber \\
b_{ik} &=& \sum_{j=1}^J(x_{ij}-c_{kj})^2\\
\text{and}\cr
\hspace{-0.5cm}\psi_{ik} &=&p_{ik},  \hspace{0.04\textwidth}\text{with  } i=1,\ldots, n:\,\, k=1,\ldots,K\cr \nonumber
\end{eqnarray}
Let us assume that $c_{kj}$ and $p_{ik}$ are known, 
replacing $(x_{ij} - c_{kj})$ with $g_{ijk}$ in \ref{obb}  we develop the following squared form:
\begin{eqnarray}
\max&\left(- \sum_{k=1}^K \sum_{i=1}^n \left( \sum_{j=1}^J g_{ijk}^2\right)p_{ik}^2\right) \nonumber \\
\mbox{s.t.}&  \sum_{k=1}^K \sum_{i=1}^n p_{ik}^2 \leq n \nonumber
\end{eqnarray}
according to the Theorem 1 we obtain:
\begin{eqnarray}
\label{mu}
\sum_{k=1}^K \sum_{i=1}^n\left( \sum_{j=1}^J g_{ijk}^2 \right)p_{ik}=\mu \sum_{k=1}^K \sum_{i=1}^n p_{ik}
\end{eqnarray}
%
%
%

%


The value of $\mu$ that optimize the \ref{mu} can be find trough the singular value decomposition of the matrix $G$, 
which is equivalent to the following Tucker3 decomposition:
\begin{eqnarray}
\label{tucker}
g_{ijk}=\sum_{r=1}^R\sum_{q=1}^Q \sum_{s=1}^S \lambda_{rqs}(u_{ir}b_{jq}v_{ks})+ e_{ijk}, \nonumber
\end{eqnarray}

with $ i=1,\ldots, n$, $j=1,\ldots, J$, $k=1,\ldots, K$.\\
Defining with: $R$ number of components of $U$, $Q$ number of components of $B$ and $S$ number of components of $V$.\\
In matrix notation:
\begin{eqnarray}
G=U\Lambda(V' \otimes B')+E
\end{eqnarray}
It can be verified that the second derivate is not positive, see section \ref{appendice:2}.\\
$\blacksquare$

The Proposition 1 and the Lemma 1 demonstrate that the Tucker3 transformation of the distance matrix $G$  minimises the JDF.
The following subsection presents an iterative algorithm to alternatively calculate $c_{kj}$ and $p_{ik}$ on one hand, and $b_{jq}$ on the other hand, until the convergence is reached. In \ref{appendice:2} we empirically demonstrate that the minimisation of the quantity  in the formula \ref{obb}  converges at least to local minima.

\subsection{Factor PD-clustering iterative algorithm}\label{dimostrazione}

Let us start considering the equation \ref{JDFhat}, where we apply the linear transformation $x_{ij}b_{jq}$ to $x_{ij}$ according to \ref{coordinate}:
\begin{eqnarray}\label{JDFstar}
\hat{JDF} = {\arg\min_{C; B}}{ \sum_{i=1}^n\sum_{q=1}^Q\sum_{k=1}^K(x_{iq}^* - c_{kq})^2 p_{ik}^2}.
\end{eqnarray}

Let us note that in formula \ref{JDFstar}:\\
$x_{ij}$ and $b_{jq}$ are the general elements of the matrices $X$ and $B$ that have been already defined in section \ref{sec:3.1};\\
$c_{kq}$ is the general element of the matrix $C$, (see eq. \ref{centri}).

It is worth to note that $C$ and $B$ are unknown matrices and $p_{ik}$ is determined as $C$ and $B$ are fixed. The problem does not
admit a direct solution and an iterative two steps procedure is required. The two alternative steps are:
\begin{itemize}
\item Linear transformation of original data;
\item PD-Clustering on transformed data.
\end{itemize}

The procedure starts with a  pseudorandomly defined centre matrix  $C$ of elements $c_{kj}$ with $k=1,\ldots , K$ and $j=1,\ldots, J$.
Then a first solution for probabilities and distance matrices is computed according to \ref{dist}.
Given the initial $C$ and $X$, the matrix $B$ is calculated; once $B$ is fixed the matrix $C$ is updated (and the values $p_{ik}$ are consequently updated). 
Last two steps are iterated until the convergence is reached: $\hat{JDF}^{(t)}-\hat{JDF}^{(t-1)} > 0$, where $t$ indicates the number of iterations.

Here under the procedure is presented according to the usual flow diagram notation:

\begin{algorithm}\begin{footnotesize}
\caption{Factor Probabilistic Distance Clustering}\label{FPDC:algo}
\begin{algorithmic}[1]
\Function{FPDC}{$X,K$}
\State $JDF \gets 1/\mathrm{eps}$           \Comment{$JDF$ is initialised to the maximum}
\State $G\gets 0$                             \Comment{Initialise the array $G$, of dimension $n\times J \times K$, to $0$}
\State $P \gets \frac{1}{k}$                             \Comment{Initialise to $\frac{1}{k}$ the probability vector $p$ of $n$ elements} 
\State $C\gets \mathrm{rand}(K,J)$ \Comment{Matrix $C_{K,J}$ is randomly initialised}

\Repeat

\For{ $k=1,K$}
\State $G_k \gets distance(X,C(k))$ \Comment{$G_k$ distances of all units from the centre $k$}
\EndFor
\State $B\gets Tucker3(G)$           \Comment{Tucker3 fun. in MatLab Toolbox N-way \cite{Che10}}
\State $X^{*} \gets XB$                  
\State $JDF0 \gets JDF$                \Comment{Current $JDF$ is stored in $JDF0$} 
\State $(C,P,JDF)\gets PDC(X^{*},K)$ \Comment{PDC() function is defined by the algorithm \ref{PDC:algo}}
\Until{$JDF0 > JDF$}

\Return{$C, P$}
\EndFunction

\end{algorithmic}\end{footnotesize}
\end{algorithm}

Remark that the Tucker3 function is in MatLab Toolbox N-way \cite{Che10}.
%
%
%
%
%
%

\section{Application on a simulated dataset}\label{applicazione}
In order to evaluate the performance of FPDC it has been applied on a simulated dataset. The dataset has been created according to Maronna and Zamar \cite{MarZam02} procedure and notations. 
\begin{figure}[h]
{\centering
\includegraphics[width=1\textwidth]{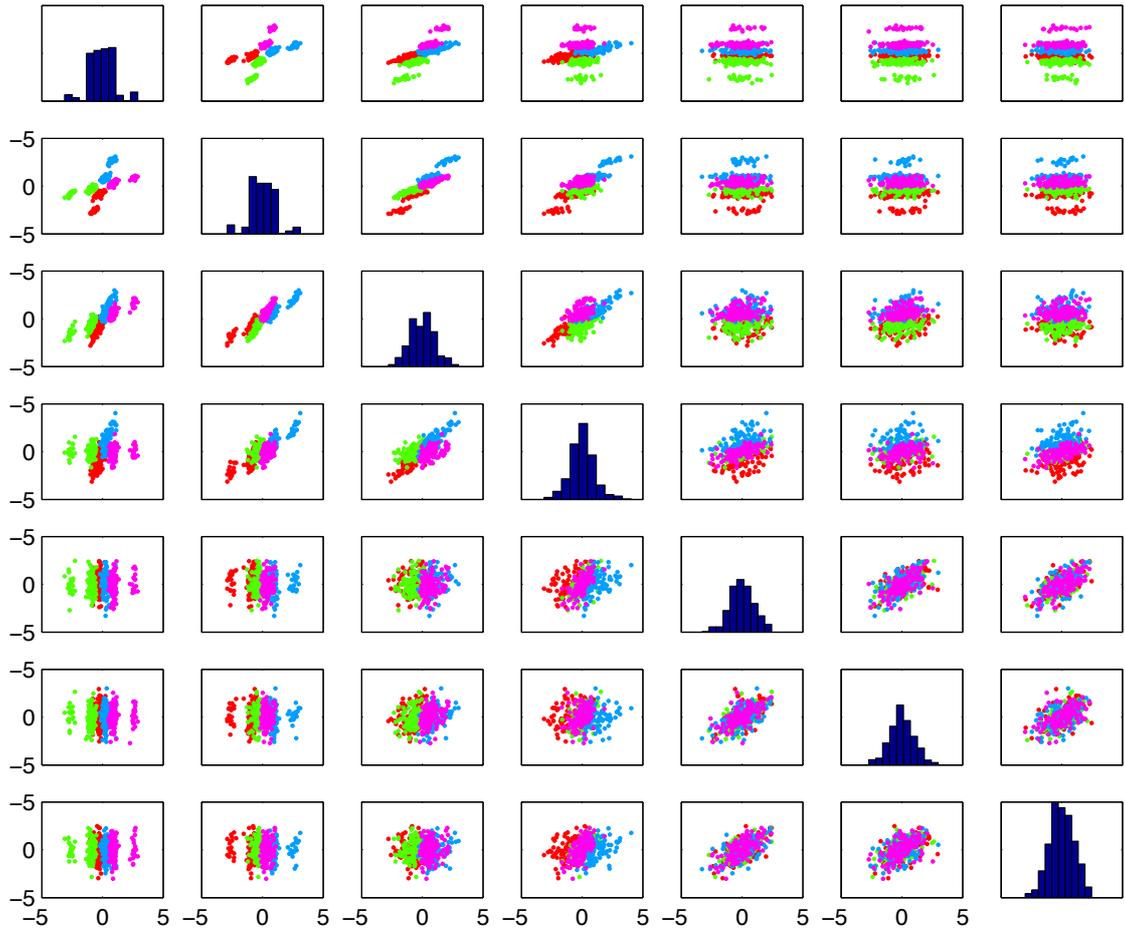}          
\caption{Scatter plot matrix of the simulated dataset. The dataset represents 4 normally generated cluster with a level of contamination of $20\%$ and correlated according to the scheme in the section \ref{applicazione}. Displayed data have been standardized. }\label{splom}}
\end{figure}

\begin{figure}[!h]
{\centering
\includegraphics[width=.75\textwidth]{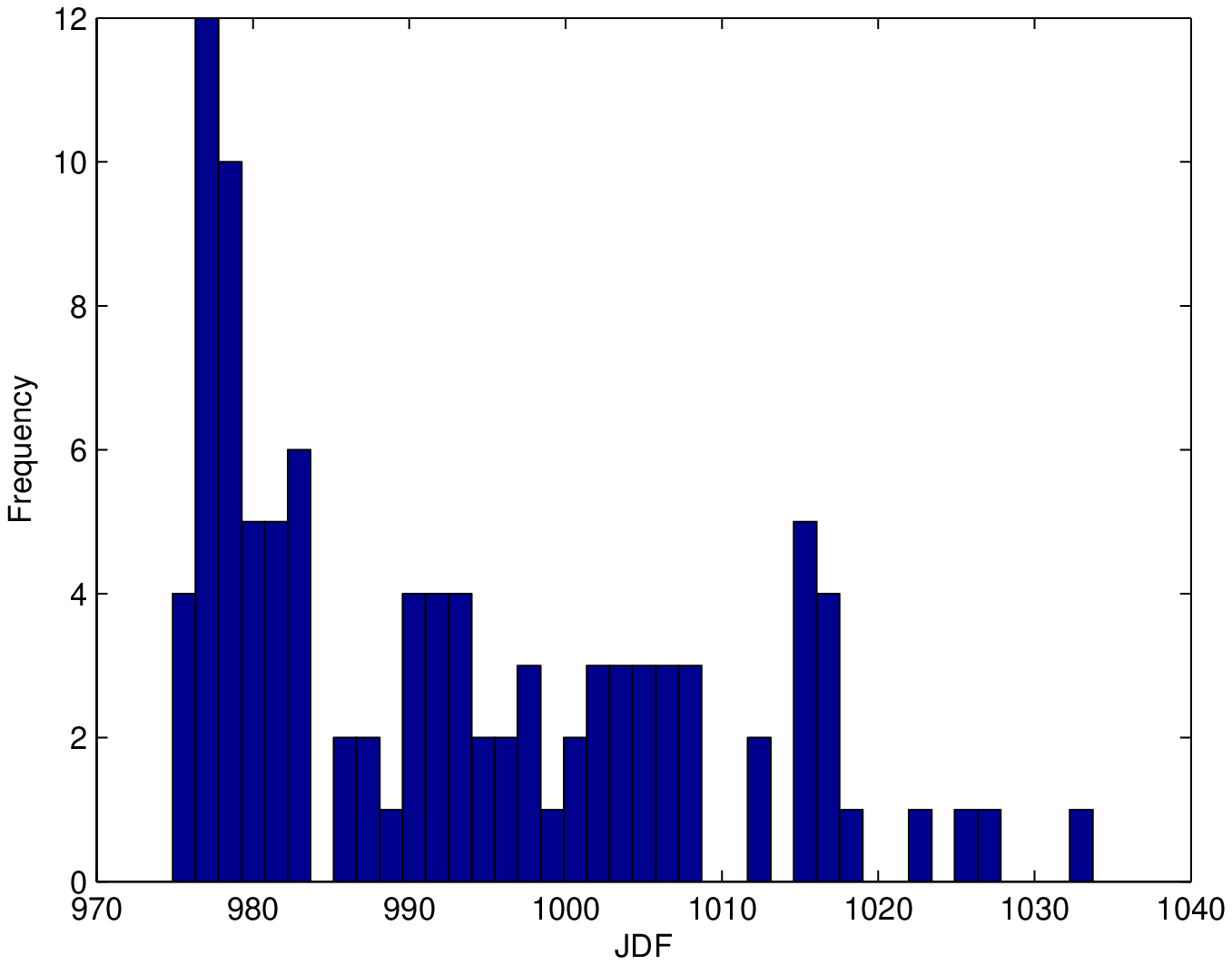}         
\caption{The bar-graph represents the distribution of JDF obtained through $100$ FPDC iterations on the simulated dataset. The picture shows  the stability of the results. The modal percentile is  $[975, 983]$ and corresponds to $59\%$ of cases. }\label{jdf} }
\end{figure}

Every cluster has been obtained generating uncorrelated normal data $x_i \sim N(0,I)$ where $I$ is a $J\times J$ identity matrix. Each element $x_i$ has been transformed into $y_i=\Sigma x_i$ where $\Sigma^2$ is a covariance matrix with $\Sigma_{jj}=1$ and $\Sigma_{jr}=\rho$ for $r \neq j$. For every cluster $100 $ vectors $y_i$ with $7$ variables  have been generated. Every cluster has been centered on points which are uniformly distributed on a hypersphere. Each cluster has been contaminated at a level $\epsilon=20\%$, cluster contamination is generated according to a normal distribution $y_i \sim N(ra_0\sqrt J, \Sigma_k)$ where $a_0$ is a unitary vector generated orthogonal to $(1,1,\ldots, 1)^T$. The parameter $r$ measures the distance between the outliers and the cluster center. To avoid that outliers overlap the elements of the clusters the minimum value of $r$ is $r_{min}=\frac{\left(1.2\sqrt{\chi^2_{J,1-\alpha}}+\sqrt{\chi^2_{J,1-\alpha}}\right)}{\sqrt J}$.  In this case we have chosen $r=4$ that verifies $r>r_{min}$.\\
In order to evaluate the stability of the results each method has been iterated $100$ times, JDF has been measured at each iteration; results are represented in fig. \ref{jdf}.

\begin{figure}[!h]
{\centering
\includegraphics[width=.75\textwidth]{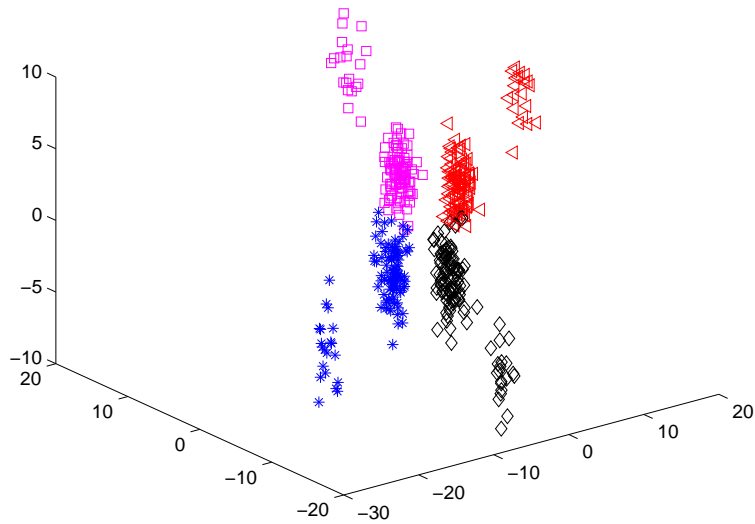}    
  \caption{ The figure shows the FPDC results of the simulated dataset composed by 4 clusters. The axes correspond to the first $3$ simulated variables see also \ref{splom}. Colors and symbols are referred to FPDC results. The misclassification
  error rate is $[0,21\%, 1,5\%]$.}\label{fpdc1}      }
\end{figure}

\begin{figure}[!h]
{\centering
\includegraphics[width=.75\linewidth]{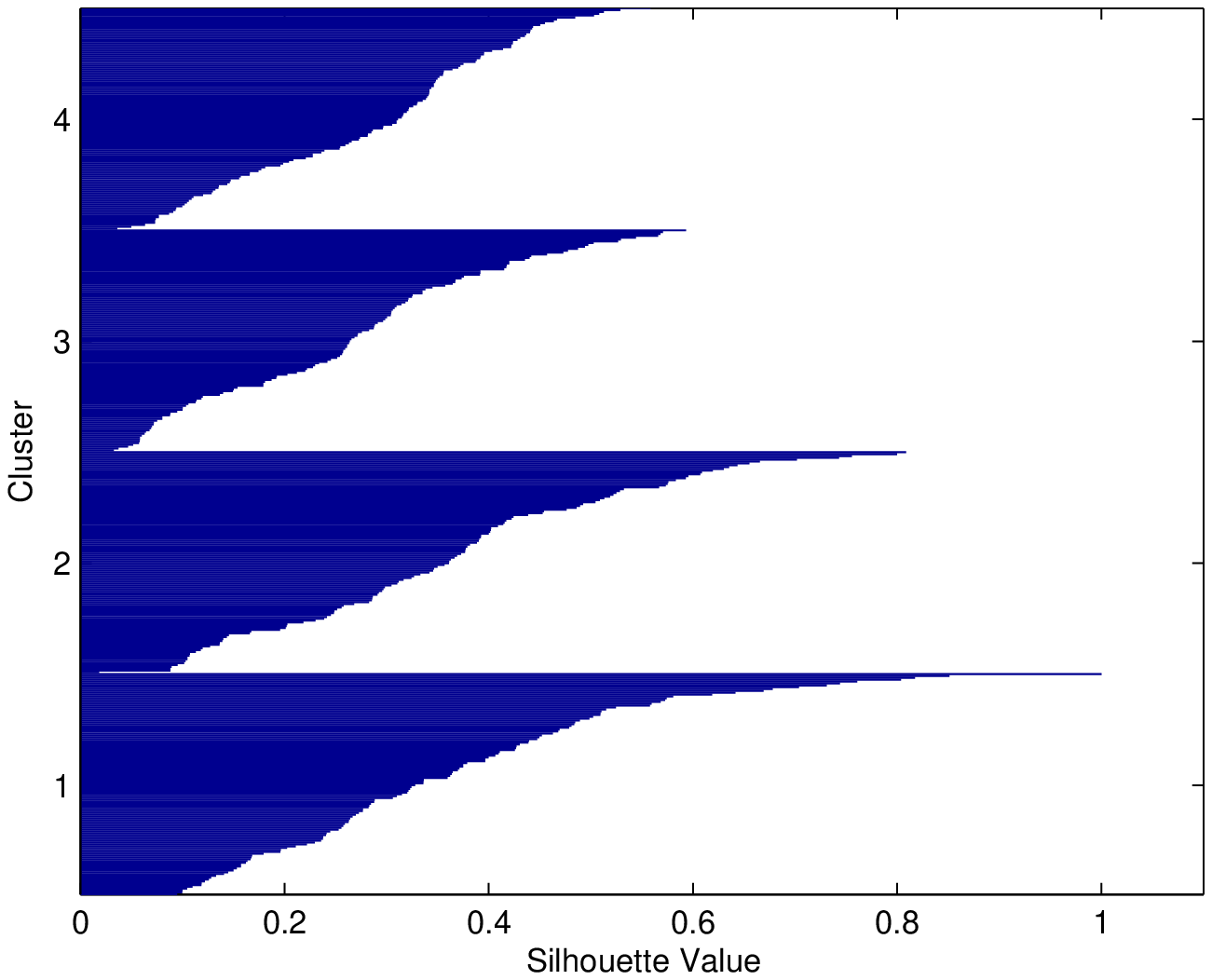}        
\caption{ The figure represents density based silhouette plot on clusters obtained in the modal value of the JDF on $100$ FPDC iterations. The graphic shows that points have been rightly classified.}\label{sil}  }
\end{figure}

The modal percentile is obtained in $59\%$ of cases, the JDF is included in the interval $[975, 983]$. In this percentile the maximum variation in clustering structure is $1\%$ that corresponds to six units. \\
In $59\%$ of cases the error term is in the interval $[0,21\%, 1,5\%]$. The clustering structure on the first three variables is represented in fig. \ref{fpdc1}.

A well known problem in cluster analysis is the validation of clustering structure. There is no index that measures clustering results because each clustering method optimizes a different function. In order to evaluate the cluster partition a density based silhouette plot (dbs) can be used. According to this method the dbs index is measured for all the observations $x_i$, all the clusters are sorted in a decreasing order with respect to dbs and plotted on a bar graph, fig. \ref{sil}. Usually euclidean distance is used to measure the distance between clusters center and each datapoint; however Euclidean distance is not suitable dealing with probabilistic clustering. A measure of dbs for probabilistic clustering method is proposed in Menardi \cite{Men11}. An adaptation of this measure for FPDC is the following one:

\begin{eqnarray}\label{dbs}
dbs_i = \frac {  log\left(\frac{p_{im_k}} {p_{im_1}}\right)}{\max_{i=1, \ldots, n} |log\left(\frac{p_{im_k}} {p_{im_1}}\right)|},
\end{eqnarray}
where $m_k$ is such that $x_i$ belongs to cluster $k$ and $m_1$ is such that $p_{im_1}$ is maximum for $m\neq m_k$.
The graphic shows that the clustering structure is correct.

 \begin{figure}[h!] 
 {\centering
 \includegraphics[width=.75\linewidth]{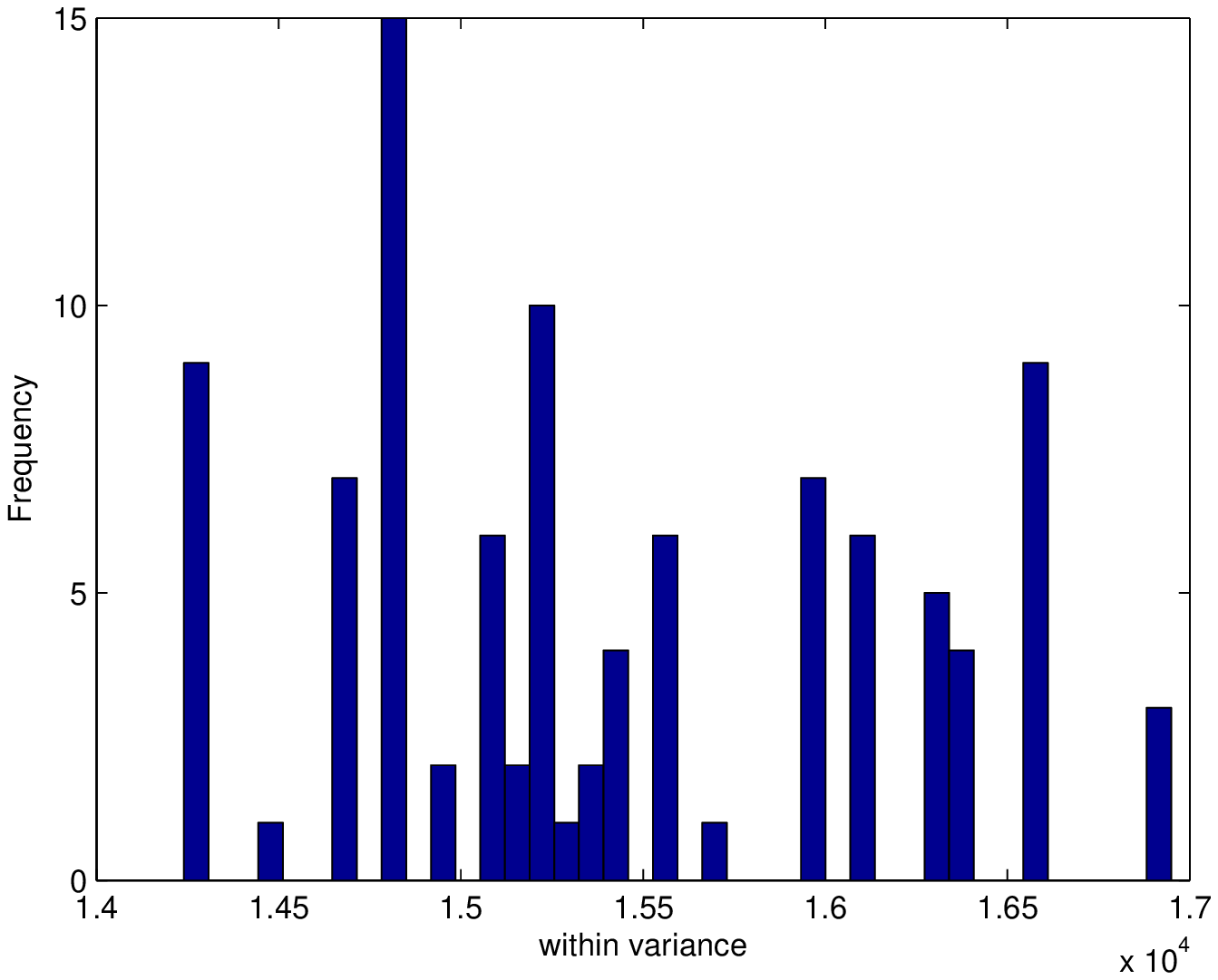}           
 \caption{  The bar-graph represents the distribution of within variance obtained through $100$ k-means iterations on the simulated dataset. The picture shows  the stability of the results. The modal percentile corresponds to $24\%$ of cases.
 }\label{varw} }
\end{figure}

Although an index that compares clustering structure does not exist, in order to point out the quality of
 FPDC the dataset has been partitioned using k-means method too. The method has been iterated $100$ times,  the within variance has been measured at each iteration, results are represented in fig. \ref{varw}.

\begin{figure}[h!]
{\centering
\includegraphics[width=.75\linewidth]{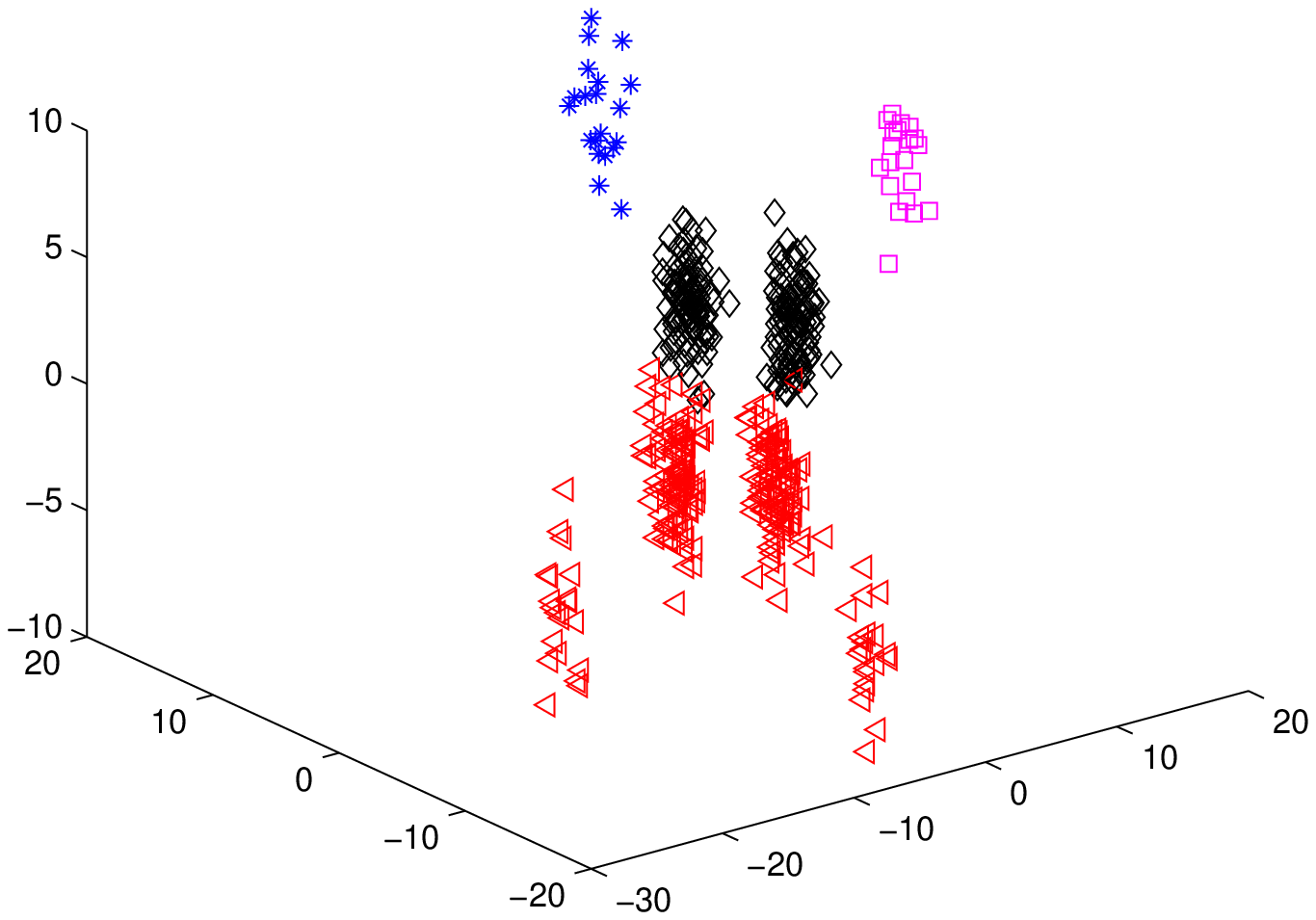}      
\caption{The figure shows the k-means results of the simulated dataset composed by 4 clusters. The axes correspond to the first $3$ simulated variables (see also \ref{splom}). Colors and symbols are referred to k-means results. The method does not find the right clustering structure in the dataset.}\label{km}    }
\end{figure}

The results have an high variability, the modal case is obtained $15\%$ of times, the first percentile is obtained $24\%$ of times. In all resulting clustering structures there is high percentage of error due to outliers.
Results obtained in the modal case are represented in fig. \ref{km}.

\section{Conclusion and perspectives}
In this paper a new factorial two-step clustering method has been brought up: Factor PD-clustering. This method can be inlaid into a new field of clustering techniques which has been developed in recent years: iterative clustering methods. Two-step clustering methods were proposed by  French school of data analysis in order to cope with some clustering issues. Thanks to computer developing, recently, iterative clustering methods have been introduced. These methods optimize a common criterion iteratively performing a linear transformation of data and a clustering optimizing a common criterion. Factor PD-clustering performs a linear transformation of data and Probabilistic D-clustering iteratively.
Probabilistic D-clustering is an iterative, distribution free, probabilistic, clustering method.
When the number of variables is large and variables are correlated PD-Clustering becomes unstable and the correlation between variables can hide the real number of clusters. A linear transformation of original variables into a reduced number of orthogonal ones using  common criteria with PD-Clustering can significantly improve the algorithm performance. 
Factor PD-clustering allows to work with large dataset improving the stability and the robustness of the method.

An important issue in the future research is the FPDC generalization to the case of categorical data. Dealing with big nominal and binary data matrices, the sparseness of data and the non-linearity in the association can be more prejudicial to the overall cluster stability. In this context, Factor clustering represents a suitable solution. Some methods have been already presented, it is worth mentioning the contributions of Hwang {\it et al.} \cite{HwaDil06} and of Palumbo and Iodice D'Enza \cite{ IodPal10}, in the case of nominal data and of binary data, respectively.

\appendix\section{Appendix}\label{app}
This appendix contains two proves: the first for  proposition 1 and the second for second order optimality condition demonstration.

\subsection{Proof of the Proposition 1} \label{App:1}
{\bf Proposition 1}
{\it Given an unknown matrix $B$ of generic element $b_{im}$ and a set of  coefficients $0 \leq  \psi_{im} \leq 1$, with $m=1,\ldots, M$ and $i=1, \ldots, n$.}
Maximising
$$-\sum_{m=1}^{M}\sum_{i=1}^nb_{im}\psi_{im}^2,$$ $\mathrm{s.t. } \sum_{m=1}^M\sum_{i=1}^n \psi^2_{im} \leq n$ 
{\it is equivalent to solve the equation}
$$\sum_{m=1}^{M}\sum_{i=1}^nb_{im}\psi_{im }= \mu\sum_{m=1}^M\sum_{i=1}^n\psi_{im},$$ 
{\it where $\mu \geq 0 $.}\\


\begin{myProof}[Proposition 1] 
To prove the proposition we introduce the Lagrangian function:
$$
\mathcal{L} =-\sum_{m=1}^{M}\sum_{i=1}^n b_{im}\psi_{im}^2 + \mu(\sum_{m=1}^{M}\sum_{i=1}^n\psi_{im}^2 -n)
$$ 
where $\mu$ is the Lagrange multiplier. Let us consider the first derivative of $\mathcal{L}$ w.r.t. $\psi_{im}$ equal to 0 :
$$
\frac{\partial\mathcal{L}}{\partial \psi_{im}} = -2\sum_{m=1}^{M}\sum_{i=1}^n b_{im}\psi_{im} + 2\mu\sum_{m=1}^{M}\sum_{i=1}^n\psi_{im} = 0 
$$
which is equivalent to
$$
\sum_{m=1}^{M}\sum_{i=1}^n{b_{im}}\psi_{im} = \mu\sum_{m=1}^{M}\sum_{i=1}^n{\psi_{im}}
$$
$ \blacksquare$
\end{myProof}

\subsection{FPD-Clustering algorithm convergence}\label{appendice:2}
\begin{figure}[h!] 
 {\centering
 \includegraphics[width=.45\linewidth]{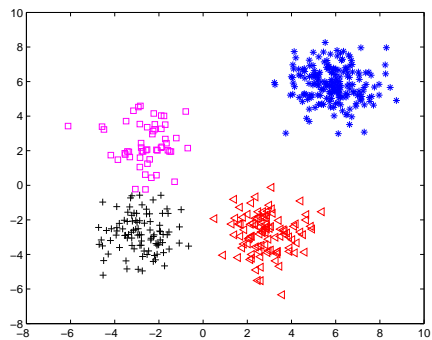}              
 \caption{  The figure represent the simulated $450 \times 2$ four clusters dataset. 
 }\label{dataset2} }
\end{figure}

In general the proof of the algorithm convergence requires the demonstration of the convexity of the objective function. Dealing with multivariate data, the analytical proof of the convexity becomes a complex issue. In most multivariate situations the empirical evidence is a satisfactory approach to verify the algorithm convergence. Moreover the high capacity of modern CPU permits to get the minimum, avoiding local minima, through the {\it multiple starts} of the algorithm.
 This section aims at empirically showing the procedure  convergence whereas a simulation study has been conducted by \cite{TorMar11}. 
The proposition states that  the convergence to a global or to a local maximum is guaranteed.
Two data sets are generated; the first one is the one used in section \ref{applicazione}. The second one is a simulated $450 \times 2$ four clusters dataset where variables are independent (see fig. \ref{dataset2}). The four clusters have been generated according to four normal distributions with different number of elements.

 \begin{figure}[h!] 
 {\centering
 \includegraphics[width=.45\linewidth]{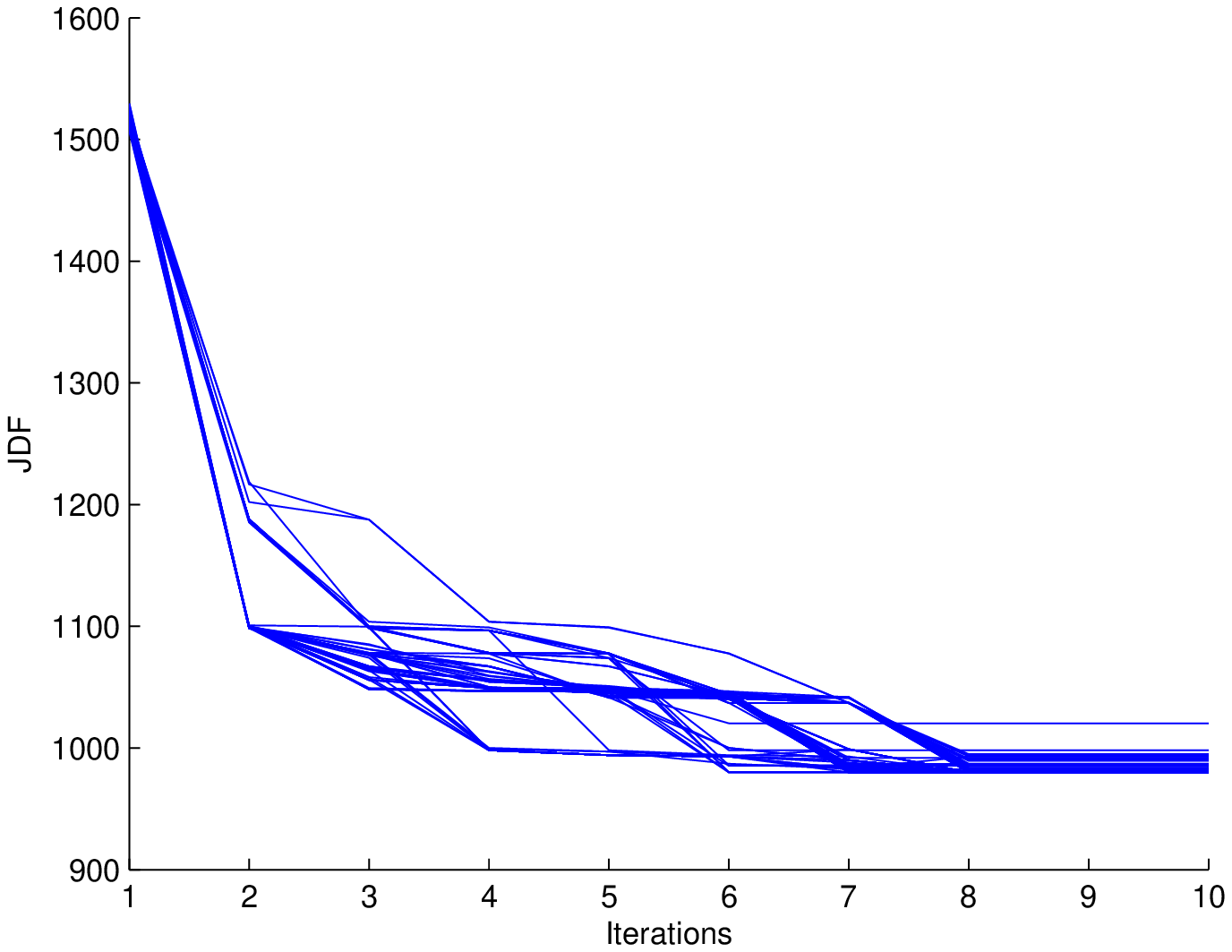}       
  \includegraphics[width=.45\linewidth]{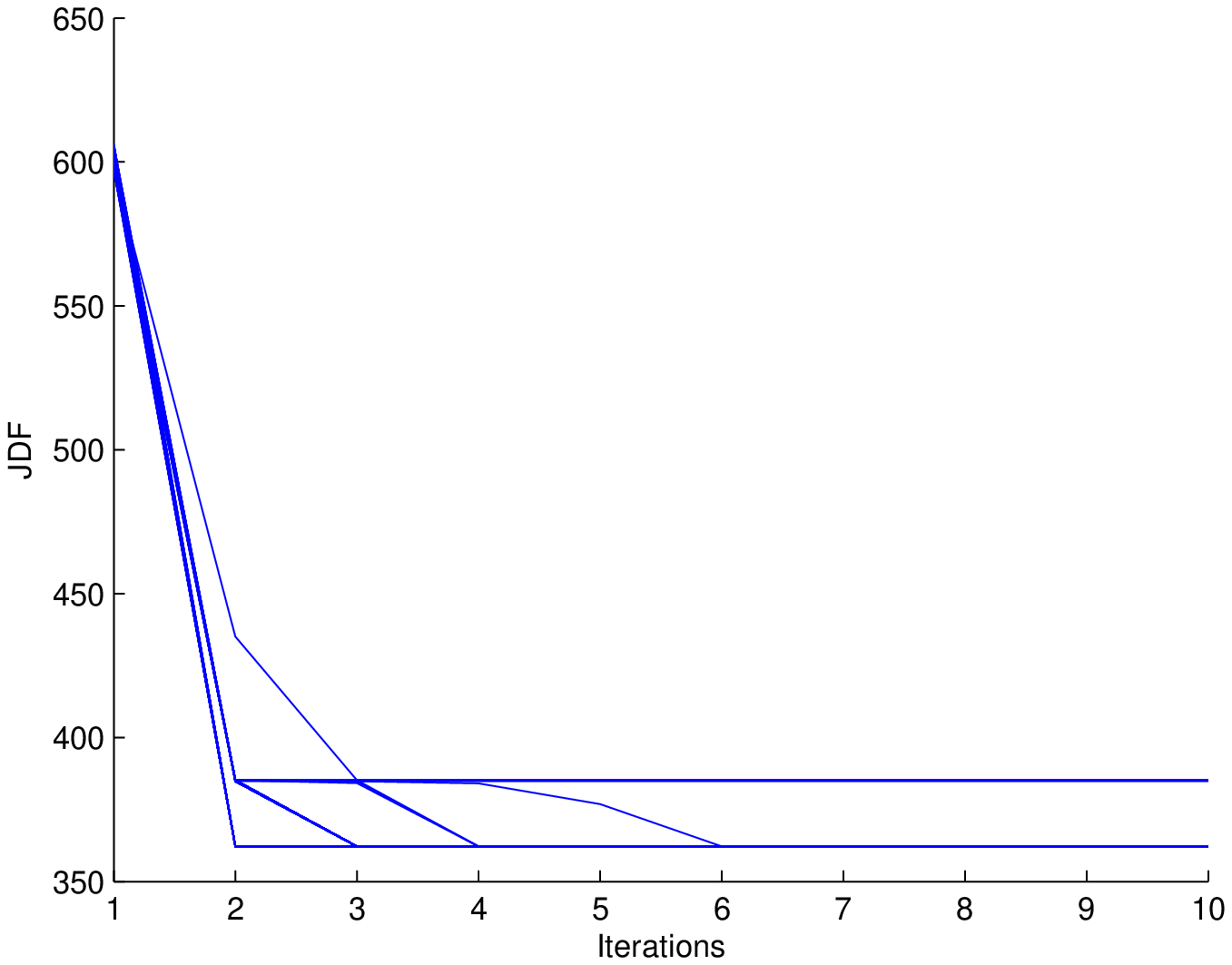}         
 \caption{  The displays represent the JDF behavior at each iteration of FPDC algorithm obtained along $100$ iterations on two simulated datasets.
 }\label{conv} }
\end{figure}
Figure \ref{conv} represents the following results: on the left-hand side the convergence of 
the  dataset one, the right-hand side of the dataset two. The
horizontal axis represents the number of iterations, the vertical refers to the
value of JDF. Each broken line represents the value of the criterion at
each iteration. When convergence is reached the line is straight and parallel to
the horizontal axis. In both cases the procedure converges in a limited number of iterations. It is worth to note that the first iteration is not counted.





\bibliographystyle{elsarticle-num}







\end{document}